\documentclass{amsart}
\usepackage{cases}
\usepackage{graphicx}
\usepackage{amssymb}
\usepackage{mathrsfs}

\newtheorem{theorem}{Theorem}[section]

\newtheorem{lemma}[theorem]{Lemma}

\numberwithin{equation}{section}

\newcommand{\bqn}{\begin{equation}}
\newcommand{\eqn}{\end{equation}}

\begin{document}
\title[]
{Note On Certain Inequalities for Neuman Means}
\author{Zai-Yin He}

\address{%
Zai-Yin He, School of Mathematics and Computation Science, Hunan City University,
Yiyang 413000, P.R. China} \

\author{Yu-Ming Chu}

\address{%
Yu-Ming Chu (Corresponding Author), School of Mathematics and Computation Science, Hunan City University,
Yiyang 413000, P.R. China}

\email{chuyuming2005@126.com}

\author{Ying-Qing Song}

\address{%
Ying-Qing Song, School of Mathematics and Computation Science, Hunan City University,
Yiyang 413000, P.R. China}

\author{Xiao-Jing Tao}

\address{%
Xiao-Jing Tao, School of Mathematics and Econometrics, Hunan University,
Changsha 410082, P.R. China}

\subjclass[2000]{26E60} \keywords{Schwab-Borchardt
mean; Neuman mean; harmonic mean; arithmetic mean; contraharmonic mean}

\begin{abstract}
In this paper, we give the explicit formulas for the Neuman means $N_{AH}$, $N_{HA}$, $N_{AC}$ and $N_{CA}$, and present the best possible upper and lower bounds for theses means in terms of the combinations of harmonic mean $H$, arithmetic mean $A$ and contraharmonic mean $C$.
\footnotesize
\end{abstract}

\maketitle

\section{Introduction}
\bigskip

Let $a, b, c\geq0$ with $ab+ac+bc\neq 0$. Then the symmetric integral $R_{F}(a, b, c)$ [1] of the first kind is defined as
\begin{equation*}
R_{F}(a,b,c)=\frac{1}{2}\int_{0}^{\infty}[(t+a)(t+b)(t+c)]^{-1/2}dt.
\end{equation*}

The degenerate case of $R_{F}$, denoted by $R_{C}$ plays an important role in the theory of special functions [1, 2], which is given by
\begin{equation*}
R_{C}(a,b)=R_{F}(a,b,b).
\end{equation*}

For $a,b>0$ with $a\neq b$, the Schwab-Borchardt mean $SB(a,b)$ [3-5] of $a$ and $b$ is given by
\begin{equation*}
SB(a,b)=\begin{cases}
\frac{\sqrt{b^2-a^2}}{\cos^{-1}{(a/b)}}, &\quad  a<b, \\
\frac{\sqrt{a^2-b^2}}{\cosh^{-1}{(a/b)}}, &\quad a>b,
\end{cases}
\end{equation*}
where $\cos^{-1}(x)$ and $\cosh^{-1}(x)=\log(x+\sqrt{x^{2}-1})$ are
the inverse cosine and inverse hyperbolic cosine functions, respectively.

Carson [6] (see also [7, (3.21)]) proved that
\begin{equation*}
SB(a,b)=\left[R_{C}\left(a^{2}, b^{2}\right)\right]^{-1}.
\end{equation*}

Recently, the Schwab-Borchardt mean has been the
subject of intensive research. In particular, many remarkable
inequalities for the Schwab-Borchardt mean and it generated means can be
found in the literature [3-5, 8-11].

Let $a>b>0$, $v=(a-b)/(a+b)\in (0, 1)$, $p\in (0, \infty)$, $q\in (0, \pi/2)$, $r\in (0, \log(2+\sqrt{3}))$ and $s\in (0, \pi/3)$ be the parameters such that $1/\cosh(p)=\cos(q)=1-v^2$, $\cosh(r)=\sec(s)=1+v^2$, $H(a,b)=2ab/(a+b)$, $G(a,b)=\sqrt{ab}$, $A(a,b)=(a+b)/2$, $Q(a,b)=\sqrt{(a^{2}+b^{2})/2}$ and $C(a,b)=(a^2+b^2)/(a+b)$ be respectively the harmonic, geometric, arithmetic, quadratic and contraharmonic means of $a$ and $b$, $S_{AH}(a,b)=SB[A(a,b), H(a,b)]$, $S_{HA}(a,b)=SB[H(a,b), A(a,b)]$, $S_{AC}(a,b)=SB[A(a,b), C(a,b)]$, $S_{CA}(a,b)=SB[C(a,b), A(a,b)]$. Then Neuman [10] gave the explicit formulas
\begin{equation}
S_{AH}(a,b)=A(a,b)\frac{\tanh(p)}{p}, \quad S_{HA}(a,b)=A(a,b)\frac{\sin q}{q},
\end{equation}
\begin{equation}
S_{CA}(a,b)=A(a,b)\frac{\sinh(r)}{r}, \quad S_{AC}(a,b)=A(a,b)\frac{\tan s}{s}.
\end{equation}

Very recently, Neuman [12] found a new mean $N(a,b)$ derived from the Schwab-Borchardt mean as follows:
\begin{equation}
N(a,b)=\frac{1}{2}\left[a+\frac{b^{2}}{SB(a,b)}\right].
\end{equation}

Let $N_{AH}(a,b)=N[A(a,b), H(a,b)]$, $N_{HA}(a,b)=N[H(a,b), A(a,b)]$, $N_{AG}(a,b)=N[A(a,b), G(a,b)]$, $N_{GA}(a,b)=N[G(a,b), A(a,b)]$,
$N_{AC}(a,b)=N[A(a,b), C(a,b)]$, $N_{CA}(a,b)=N[C(a,b), A(a,b)]$, $N_{AQ}(a,b)=N[A(a,b), Q(a,b)]$ and $N_{QA}(a,b)=N[Q(a,b), A(a,b)]$ be the Neuman means. Then Neuman [12] proved that
\begin{equation*}
G(a,b)<N_{AG}(a,b)<N_{GA}(a,b)<A(a,b)<N_{QA}(a,b)<N_{AQ}(a,b)<Q(a,b)
\end{equation*}
for all $a, b>0$ with $a\neq b$, and the double inequalities
$$\alpha_{1} A(a,b)+(1-\alpha_{1})G(a,b)<N_{GA}(a,b)<\beta_{1} A(a,b)+(1-\beta_{1})G(a,b),$$
$$\alpha_{2} Q(a,b)+(1-\alpha_{2})A(a,b)<N_{AQ}(a,b)<\beta_{2} Q(a,b)+(1-\beta_{2})A(a,b),$$
$$\alpha_{3} A(a,b)+(1-\alpha_{3})G(a,b)<N_{AG}(a,b)<\beta_{3} A(a,b)+(1-\beta_{3})G(a,b),$$
$$\alpha_{4} Q(a,b)+(1-\alpha_{4})A(a,b)<N_{QA}(a,b)<\beta_{4} Q(a,b)+(1-\beta_{4})A(a,b)$$
hold for all $a, b>0$ with $a\neq b$ if and only if $\alpha_{1}\leq 2/3$, $\beta_{1}\geq \pi/4$,
$\alpha_{2}\leq 2/3$, $\beta_{2}\geq (\pi-2)/[4(\sqrt{2}-1)]=0.689\ldots$, $\alpha_{3}\leq 1/3$, $\beta_{3}\geq 1/2$,
$\alpha_{4}\leq 1/3$ and $\beta_{4}\geq [\log(1+\sqrt{2})+\sqrt{2}-2]/[2(\sqrt{2}-1)]=0.356\ldots$.

Zhang et. al. [13] presented the best possible parameters $\alpha_{1}, \alpha_{2}, \beta_{1}, \beta_{2}\in [0, 1/2]$ and
$\alpha_{3}, \alpha_{4}, \beta_{3}, \beta_{4}\in [1/2, 1]$ such that the double inequalities
$$G(\alpha_{1}a+(1-\alpha_{1})b, \alpha_{1}b+(1-\alpha_{1})a)<N_{AG}(a,b)<G(\beta_{1}a+(1-\beta_{1})b, \beta_{1}b+(1-\beta_{1})a),$$
$$G(\alpha_{2}a+(1-\alpha_{2})b, \alpha_{2}b+(1-\alpha_{2})a)<N_{GA}(a,b)<G(\beta_{2}a+(1-\beta_{2})b, \beta_{2}b+(1-\beta_{2})a),$$
$$Q(\alpha_{3}a+(1-\alpha_{3})b, \alpha_{3}b+(1-\alpha_{3})a)<N_{QA}(a,b)<Q(\beta_{3}a+(1-\beta_{3})b, \beta_{3}b+(1-\beta_{3})a),$$
$$Q(\alpha_{4}a+(1-\alpha_{4})b, \alpha_{4}b+(1-\alpha_{4})a)<N_{AQ}(a,b)<Q(\beta_{4}a+(1-\beta_{4})b, \beta_{4}b+(1-\beta_{4})a)$$
hold for all $a, b>0$ with $a\neq b$.

In [14], the authors found the greatest values $\alpha_{1}$, $\alpha_{2}$, $\alpha_{3}$, $\alpha_{4}$, $\alpha_{5}$, $\alpha_{6}$, $\alpha_{7}$,
$\alpha_{8}$ and the least values $\beta_{1}$, $\beta_{2}$, $\beta_{3}$, $\beta_{4}$, $\beta_{5}$, $\beta_{6}$, $\beta_{7}$, $\beta_{8}$ such that the double inequalities
$$A^{\alpha_{1}}(a,b)G^{1-\alpha_{1}}(a,b)<N_{GA}(a,b)<A^{\beta_{1}}(a,b)G^{1-\beta_{1}}(a,b),$$
$$\frac{\alpha_{2}}{G(a,b)}+\frac{1-\alpha_{2}}{A(a,b)}<\frac{1}{N_{GA}(a,b)}<\frac{\beta_{2}}{G(a,b)}+\frac{1-\beta_{2}}{A(a,b)},$$
$$A^{\alpha_{3}}(a,b)G^{1-\alpha_{3}}(a,b)<N_{AG}(a,b)<A^{\beta_{3}}(a,b)G^{1-\beta_{3}}(a,b),$$
$$\frac{\alpha_{4}}{G(a,b)}+\frac{1-\alpha_{4}}{A(a,b)}<\frac{1}{N_{AG}(a,b)}<\frac{\beta_{4}}{G(a,b)}+\frac{1-\beta_{4}}{A(a,b)},$$
$$Q^{\alpha_{5}}(a,b)A^{1-\alpha_{5}}(a,b)<N_{AQ}(a,b)<Q^{\beta_{5}}(a,b)A^{1-\beta_{5}}(a,b),$$
$$\frac{\alpha_{6}}{A(a,b)}+\frac{1-\alpha_{6}}{Q(a,b)}<\frac{1}{N_{AQ}(a,b)}<\frac{\beta_{6}}{A(a,b)}+\frac{1-\beta_{6}}{Q(a,b)},$$
$$Q^{\alpha_{7}}(a,b)A^{1-\alpha_{7}}(a,b)<N_{QA}(a,b)<Q^{\beta_{7}}(a,b)A^{1-\beta_{7}}(a,b),$$
$$\frac{\alpha_{8}}{A(a,b)}+\frac{1-\alpha_{8}}{Q(a,b)}<\frac{1}{N_{QA}(a,b)}<\frac{\beta_{8}}{A(a,b)}+\frac{1-\beta_{8}}{Q(a,b)}$$
hold for all $a, b>0$ with $a\neq b$.

The main purpose of this paper is to give the explicit formulas for the Neuman means $N_{AH}$, $N_{HA}$, $N_{AC}$ and $N_{CA}$, and present the best possible upper and lower bounds for theses means in terms of the combinations of harmonic, arithmetic and contraharmonic means. Our main results are the following Theorems 1.1-1.3.

\medskip
\begin{theorem}  Let $a>b>0$, $v=(a-b)/(a+b)\in (0, 1)$, $p\in (0, \infty)$, $q\in (0, \pi/2)$, $r\in (0, \log(2+\sqrt{3}))$ and $s\in (0, \pi/3)$ be the parameters such that $1/\cosh(p)=\cos(q)=1-v^2$, $\cosh(r)=\sec(s)=1+v^2$. Then we have
\begin{equation}
N_{AH}(a,b)=\frac{1}{2}A(a,b)\left[1+\frac{2p}{\sinh(2p)}\right],
\end{equation}
\begin{equation}
N_{HA}(a,b)=\frac{1}{2}A(a,b)\left[\cos(q)+\frac{q}{\sin(q)}\right],
\end{equation}
\begin{equation}
N_{CA}(a,b)=\frac{1}{2}A(a,b)\left[\cosh(r)+\frac{r}{\sinh(r)}\right],
\end{equation}
\begin{equation}
N_{AC}(a,b)=\frac{1}{2}A(a,b)\left[1+\frac{2s}{\sin(2s)}\right]
\end{equation}
and
\begin{equation}
H(a,b)<N_{AH}(a,b)<N_{HA}(a,b)<A(a,b)
\end{equation}
\begin{equation*}
<N_{CA}(a,b)<N_{AC}(a,b)<C(a,b).
\end{equation*}
\end{theorem}

\medskip
\begin{theorem} The double inequalities
\begin{equation}
\alpha_{1}A(a,b)+(1-\alpha_{1})H(a,b)<N_{AH}(a,b)<\beta_{1}A(a,b)+(1-\beta_{1})H(a,b),
\end{equation}
\begin{equation}
\alpha_{2}A(a,b)+(1-\alpha_{2})H(a,b)<N_{HA}(a,b)<\beta_{2}A(a,b)+(1-\beta_{2})H(a,b),
\end{equation}
\begin{equation}
\alpha_{3}C(a,b)+(1-\alpha_{3})A(a,b)<N_{CA}(a,b)<\beta_{3}C(a,b)+(1-\beta_{3})A(a,b),
\end{equation}
\begin{equation}
\alpha_{4}C(a,b)+(1-\alpha_{4})A(a,b)<N_{AC}(a,b)<\beta_{4}C(a,b)+(1-\beta_{4})A(a,b)
\end{equation}
hold for all $a, b>0$ with $a\neq b$ if and only if $\alpha_{1}\leq
1/3$, $\beta_{1}\geq 1/2$, $\alpha_{2}\leq
2/3$, $\beta_{2}\geq \pi/4=0.7853\ldots$, $\alpha_{3}\leq
1/3$, $\beta_{3}\geq\sqrt{3}\log(2+\sqrt{3})/6=0.3801\ldots$,
$\alpha_{4}\leq 2/3$ and $\beta_{4}\geq(4\sqrt{3}\pi-9)/18=0.7901\ldots$.
\end{theorem}

\medskip
\begin{theorem} The double inequalities
\begin{equation}
\frac{\alpha_{5}}{H(a,b)}+\frac{1-\alpha_{5}}{A(a,b)}<\frac{1}{N_{AH}(a,b)}<\frac{\beta_{5}}{H(a,b)}+\frac{1-\beta_{5}}{A(a,b)},
\end{equation}
\begin{equation}
\frac{\alpha_{6}}{H(a,b)}+\frac{1-\alpha_{6}}{A(a,b)}<\frac{1}{N_{HA}(a,b)}<\frac{\beta_{6}}{H(a,b)}+\frac{1-\beta_{6}}{A(a,b)},
\end{equation}
\begin{equation}
\frac{\alpha_{7}}{A(a,b)}+\frac{1-\alpha_{7}}{C(a,b)}<\frac{1}{N_{CA}(a,b)}<\frac{\beta_{7}}{A(a,b)}+\frac{1-\beta_{7}}{C(a,b)},
\end{equation}
\begin{equation}
\frac{\alpha_{8}}{A(a,b)}+\frac{1-\alpha_{8}}{C(a,b)}<\frac{1}{N_{AC}(a,b)}<\frac{\beta_{8}}{A(a,b)}+\frac{1-\beta_{8}}{C(a,b)},
\end{equation}
hold for all $a, b>0$ with $a\neq b$ if and only if $\alpha_{5}\leq
0$, $\beta_{5}\geq 2/3$, $\alpha_{6}\leq 0$,
$\beta_{6}\geq 1/3$, $\alpha_{7}\leq
[2\sqrt{3}-\log(2+\sqrt{3})]/[2\sqrt{3}+\log(2+\sqrt{3})]=0.4490\ldots$,
$\beta_{7}\geq 2/3$, $\alpha_{8}\leq
(9\sqrt{3}-4\pi)/(3\sqrt{3}+4\pi)=0.1701\ldots$ and $\beta_{8}\geq 1/3$.
\end{theorem}

\bigskip
\section {Lemmas}
\bigskip

In order to prove our main results we need several lemmas, which we present in this section.

\setcounter{equation}{0}

\medskip
\begin{lemma} (See [15, Theorem 1.25]) For $-\infty<a<b<\infty$,
let $f,g:[a,b]\rightarrow{\mathbb{R}}$ be continuous on $[a,b]$, and
be differentiable on $(a,b)$, let $g'(x)\neq 0$ on $(a,b)$. If
$f^{\prime}(x)/g^{\prime}(x)$ is increasing (decreasing) on $(a,b)$,
then so are
$$\frac{f(x)-f(a)}{g(x)-g(a)}\ \ \mbox{and}\ \ \frac{f(x)-f(b)}{g(x)-g(b)}.$$
If $f^{\prime}(x)/g^{\prime}(x)$ is strictly monotone, then the
monotonicity in the conclusion is also strict.
\end{lemma}

\medskip
\begin{lemma}(See [16, Lemma 1.1]) Suppose that the power
series $f(x)={\sum^{\infty}_{n=0}}a_{n}x^{n}$ and
$g(x)={\sum^{\infty}_{n=0}}b_{n}x^{n}$ have the radius of convergence $r>0$ and $a_{n}, b_{n}>0$ for all $n\geq 0$.
If the sequence $\{a_{n}/b_{n}\}$ is (strictly) increasing (decreasing) for all $n\geq 0$, then the function $f(x)/g(x)$ is also
(strictly) increasing (decreasing)  on $(0, r)$.
\end{lemma}

\medskip
\begin{lemma} (See [12, Theorem 4.1]) If $a>b$, then
\begin{equation*}
N(b,a)>N(a,b).
\end{equation*}
\end{lemma}

\medskip
\begin{lemma} The function
\begin{equation*}
\varphi_{1}(t)=\frac{\sinh(2t)-4\sinh(t)+2t}{\sinh(2t)-2\sinh(t)}
\end{equation*}
is strictly increasing from $(0, \infty)$ onto $(2/3, 1)$.
\end{lemma}
{\em Proof.} Making use of power series expansion we get
\begin{equation}
\varphi_{1}(t)=\frac{\sum_{n=1}^{\infty}\frac{2^{2n+1}-4}{(2n+1)!}t^{2n+1}}{\sum_{n=1}^{\infty}\frac{2^{2n+1}-2}{(2n+1)!}t^{2n+1}}
=\frac{\sum_{n=0}^{\infty}\frac{2^{2n+3}-4}{(2n+3)!}t^{2n}}{\sum_{n=0}^{\infty}\frac{2^{2n+3}-2}{(2n+3)!}t^{2n}}.
\end{equation}

Let
\begin{equation}
a_{n}=\frac{2^{2n+3}-4}{(2n+3)!}, \quad b_{n}=\frac{2^{2n+3}-2}{(2n+3)!}.
\end{equation}

Then
\begin{equation}
a_{n}>0, \quad b_{n}>0
\end{equation}
and $a_{n}/b_{n}=1-1/(2^{2n+2}-1)$ is strictly increasing for all $n\geq 0$.

Note that
\begin{equation}
\varphi_{1}(0^{+})=\frac{a_{0}}{b_{0}}=\frac{2}{3}, \quad \varphi_{1}(\infty)=\lim_{n\rightarrow \infty}\frac{a_{n}}{b_{n}}=1.
\end{equation}

Therefore, Lemma 2.4 follows easily from Lemma 2.2 and (2.1)-(2.4) together with the monotonicity of the sequence $\{a_{n}/b_{n}\}$.  $\Box$

\medskip
\begin{lemma} The function
\begin{equation*}
\varphi_{2}(t)=\frac{2t-\sin(2t)}{\sin t(1-\cos t)}
\end{equation*}
is strictly increasing from $(0, \pi/2)$ onto $(8/3, \pi)$.
\end{lemma}
{\em Proof.} Let $f_{1}(t)=2t-\sin(2t)$ and $g_{1}(t)=\sin t(1-\cos t)$. Then simple computations lead to
\begin{equation}
\varphi_{2}(t)=\frac{f_{1}(t)-f_{1}(0)}{g_{1}(t)-g_{1}(0)}
\end{equation}
and $f'_{1}(t)/g'_{1}(t)=4[1-1/(2+1/\cos t)]$ is strictly increasing on $(0, \pi/2)$.

Note that
\begin{equation}
\varphi_{2}(0^{+})=\lim_{t\rightarrow 0^{+}}\frac{f'_{1}(t)}{g'_{1}(t)}=\frac{8}{3}, \quad \varphi_{2}(\pi/2)=\pi.
\end{equation}

Therefore, Lemma 2.5 follows from Lemma 2.1, (2,5), (2.6) and the monotonicity of $f'_{1}(t)/g'_{1}(t)$. $\Box$

\medskip
\begin{lemma} The function
\begin{equation*}
\varphi_{3}(t)=\frac{\sinh(t)\cosh(t)-t}{[\sinh(t)\cosh(t)+t](\cosh(t)-1)}
\end{equation*}
is strictly decreasing from $(0, \infty)$ onto $(0, 2/3)$.
\end{lemma}
{\em Proof.} Simple computations lead to
\begin{equation}
\varphi_{3}(t)=\frac{2\sinh(2t)-4t}{\sinh(3t)+4t\cosh(t)+\sinh(t)-2\sinh(2t)-4t}
\end{equation}
\begin{equation*}
=\frac{\sum_{n=0}^{\infty}\frac{2^{2n+4}}{(2n+3)!}t^{2n}}{\sum_{n=0}^{\infty}\frac{3^{2n+3}-2^{2n+4}+8n+13}{(2n+3)!}t^{2n}}.
\end{equation*}

Let
\begin{equation}
a_{n}=\frac{2^{2n+4}}{(2n+3)!}, \quad b_{n}=\frac{3^{2n+3}-2^{2n+4}+8n+13}{(2n+3)!}.
\end{equation}

Then
\begin{equation}
a_{n}>0, \quad b_{n}>0
\end{equation}
and
\begin{equation}
\frac{a_{n+1}}{b_{n+1}}-\frac{a_{n}}{b_{n}}=-\frac{2^{2n+4}\left(5\times 3^{2n+3}-24n-31\right)}{\left(3^{2n+5}-2^{2n+6}+8n+21\right)\left(3^{2n+3}-2^{2n+4}+8n+13\right)}<0
\end{equation}
for all $n\geq 0$.

Note that
\begin{equation}
\varphi_{3}(0^{+})=\frac{a_{0}}{b_{0}}=\frac{2}{3}, \quad \varphi_{3}(\infty)=\lim_{n\rightarrow \infty}\frac{a_{n}}{b_{n}}=0.
\end{equation}

Therefore, Lemma 2.6 follows easily from (2.7)-(2.11) and Lemma 2.2.  $\Box$

\medskip
\begin{lemma} The function
\begin{equation*}
f(t)=9\cos t+\frac{t}{\sin t}
\end{equation*}
is strictly decreasing on the interval $(0, \pi/2)$.
\end{lemma}
{\em Proof.} Let $f_{2}(t)=9\sin t\cos t+t$ and $g_{2}(t)=\sin t$. Then simple computations lead to
\begin{equation}
f(t)=\frac{f_{2}(t)-f_{2}(0)}{g_{2}(t)-g_{2}(0)},
\end{equation}
\begin{equation*}
\frac{f'_{2}(t)}{g'_{2}(t)}=\frac{18\cos^{2}t-8}{\cos t}
\end{equation*}
and
\begin{equation}
\left[\frac{f'_{2}(t)}{g'_{2}(t)}\right]'=-\frac{2\sin t(9\cos^{2}+4)}{\cos^{2}(t)}<0
\end{equation}
for $t\in (0, \pi/2)$.

Therefore, Lemma 2.7 follows easily from (2.12) and (2.13) together with Lemma 2.1.  $\Box$

\medskip
\begin{lemma} The function
\begin{equation*}
\varphi_{4}(t)=\frac{\sin t\cos t-t}{(t+\sin t\cos t)(1-\cos t)}
\end{equation*}
is strictly decreasing from $(0, \pi/2)$ onto $(-1, -2/3)$.
\end{lemma}
{\em Proof.} Let $f_{3}(t)=\sin t\cos t-t$ and $g_{3}(t)=(t+\sin t\cos t)(1-\cos t)$. Then simple computations lead to
\begin{equation}
\varphi_{4}(t)=\frac{f_{3}(t)}{g_{3}(t)}=\frac{f_{3}(t)-f_{3}(0)}{g_{3}(t)-g_{3}(0)},
\end{equation}
\begin{equation}
\frac{f'_{3}(t)}{g'_{3}(t)}=\frac{f'_{3}(t)-f'_{3}(0)}{g'_{3}(t)-g'_{3}(0)}
\end{equation}
and
\begin{equation}
\frac{f''_{3}(t)}{g''_{3}(t)}=\frac{4}{4-\left(9\cos t+\frac{t}{\sin t}\right)}.
\end{equation}

Note that
\begin{equation}
\varphi_{4}(0^{+})=\lim_{t\rightarrow 0^{+}}\frac{f''_{3}(t)}{g''_{3}(t)}=-\frac{2}{3}, \quad \varphi_{4}\left(\frac{\pi}{2}\right)=-1.
\end{equation}

Therefore, Lemma 2.8 follows from Lemma 2.1 and Lemma 2.7 together with (2.14)-(2.17). $\Box$

\bigskip
\section{Proofs of Theorems 1.1-1.3}
\bigskip

{\bf Proof of Theorem 1.1}. It follows from (1.1)-(1.3) we clearly see that
\begin{equation*}
N_{AH}(a,b)=\frac{1}{2}\left[A(a,b)+\frac{H^2(a,b)}{S_{AH}(a,b)}\right]=\frac{1}{2}A(a,b)\left[1+(1-v^{2})^{2}\frac{p}{\tanh(p)}\right]
\end{equation*}
\begin{equation*}
=\frac{1}{2}A(a,b)\left[1+\frac{p}{\tanh(p)\cosh^{2}(p)}\right]=\frac{1}{2}A(a,b)\left[1+\frac{2p}{\sinh(2p)}\right],
\end{equation*}
\begin{equation*}
N_{HA}(a,b)=\frac{1}{2}\left[H(a,b)+\frac{A^2(a,b)}{S_{HA}(a,b)}\right]
\end{equation*}
\begin{equation*}
=\frac{1}{2}A(a,b)\left[(1-v^{2})+\frac{q}{\sin q}\right]=\frac{1}{2}A(a,b)\left[\cos q+\frac{q}{\sin q}\right],
\end{equation*}
\begin{equation*}
N_{CA}(a,b)=\frac{1}{2}\left[C(a,b)+\frac{A^2(a,b)}{S_{CA}(a,b)}\right]
\end{equation*}
\begin{equation*}
=\frac{1}{2}A(a,b)\left[(1+v^{2})+\frac{r}{\sinh(r)}\right]=\frac{1}{2}A(a,b)\left[\cosh(r)+\frac{r}{\sinh(r)}\right],
\end{equation*}
\begin{equation*}
N_{AC}(a,b)=\frac{1}{2}\left[A(a,b)+\frac{C^2(a,b)}{S_{AC}(a,b)}\right]=\frac{1}{2}A(a,b)\left[1+(1+v^{2})^{2}\frac{s}{\tan(s)}\right]
\end{equation*}
\begin{equation*}
=\frac{1}{2}A(a,b)\left[1+\frac{s}{\tan(s)\cos^{2}s}\right]=\frac{1}{2}A(a,b)\left[1+\frac{2s}{\sin(2s)}\right].
\end{equation*}

Inequalities (1.8) follows easily from $H(a,b)<A(a,b)<C(a,b)$ and Lemma 2.3 together with the fact that $N_{KL}(a,b)$ is a mean of $K(a,b)$ and $L(a,b)$ for $K(a,b), L(a,b)\in \{H(a,b), A(a,b), C(a,b)\}$.  $\Box$

\bigskip
\bigskip
{\bf Proof of Theorem 1.2}.  Without loss of generality, we assume that $a>b$. Let $v=(a-b)/(a+b)\in (0, 1)$, $p\in (0, \infty)$, $q\in (0, \pi/2)$, $r\in (0, \log(2+\sqrt{3}))$ and $s\in (0, \pi/3)$ be the parameters such that $1/\cosh(p)=\cos(q)=1-v^2$, $\cosh(r)=\sec(s)=1+v^2$. Then from (1.4)- (1.7) we have
\begin{equation}
\frac{N_{AH}(a,b)-H(a,b)}{A(a,b)-H(a,b)}=\frac{\frac{1}{2}\left[1+\frac{2p}{\sinh(2p)}\right]-(1-v^{2})}{v^{2}}
\end{equation}
\begin{equation*}
=\frac{\frac{1}{2}\left[1+\frac{2p}{\sinh(2p)}\right]-\frac{1}{\cosh(p)}}{1-\frac{1}{\cosh(p)}}=\varphi_{1}(p),
\end{equation*}
\begin{equation}
\frac{N_{HA}(a,b)-H(a,b)}{A(a,b)-H(a,b)}=\frac{\frac{1}{2}\left[\cos q+\frac{q}{\sin q}\right]-(1-v^{2})}{v^{2}}
\end{equation}
\begin{equation*}
=\frac{\frac{1}{2}\left[\cos q+\frac{q}{\sin q}\right]-\cos q}{1-\cos q}=\frac{1}{4}\varphi_{2}(q),
\end{equation*}
\begin{equation}
\frac{N_{CA}(a,b)-A(a,b)}{C(a,b)-A(a,b)}=\frac{\frac{1}{2}\left[\cosh(r)+\frac{r}{\sinh(r)}\right]-1}{v^{2}}
\end{equation}
\begin{equation*}
=\frac{\frac{1}{2}\left[\cosh(r)+\frac{r}{\sinh(r)}\right]-1}{\cosh(r)-1}=\frac{1}{2}\varphi_{1}(r).
\end{equation*}

\begin{equation}
\frac{N_{AC}(a,b)-A(a,b)}{C(a,b)-A(a,b)}=\frac{\frac{1}{2}\left[1+\frac{2s}{\sin(2s)}\right]-1}{v^{2}}
\end{equation}
\begin{equation*}
=\frac{\frac{1}{2}\left[1+\frac{2s}{\sin(2s)}\right]-1}{\sec(s)-1}=\frac{1}{4}\varphi_{2}(s),
\end{equation*}
where the functions $\varphi_{1}$ and $\varphi_{2}$ are defined as in Lemmas 2.4 and 2.5, respectively.

Note that
\begin{equation}
\varphi_{1}[\log(2+\sqrt{3})]=\sqrt{3}\log(2+\sqrt{3})/6
\end{equation}
and
\begin{equation}
\varphi_{2}\left(\frac{\pi}{3}\right)=\frac{8\sqrt{3}\pi-18}{9}.
\end{equation}

Therefore, inequality (1.9) holds for all $a, b>0$ with $a\neq b$ if and only if $\alpha_{1}\leq 1/3$ and $\beta_{1}\geq 1/2$ follows from (3.1) and Lemma 2.4, inequality (1.10) holds for all $a, b>0$ with $a\neq b$ if and only if
 $\alpha_{2}\leq 2/3$ and $\beta_{2}\geq \pi/4$ follows from (3.2) and Lemma 2.5, inequality (1.11) holds for all $a, b>0$ with $a\neq b$ if and only if $\alpha_{3}\leq 1/3$ and $\beta_{3}\geq\sqrt{3}\log(2+\sqrt{3})/6$ follows from (3.3) and (3.5) together with Lemma 2.4, and inequality (1.12) holds for all $a, b>0$ with $a\neq b$ if and only if $\alpha_{4}\leq 2/3$ and $\beta_{4}\geq(4\sqrt{3}\pi-9)/18$ follows from (3.4) and (3.6) together with Lemma 2.5.  $\Box$

\bigskip
\bigskip
{\bf Proof of Theorem 1.3}.  Without loss of generality, we assume that $a>b$. Let $v=(a-b)/(a+b)\in (0, 1)$, $p\in (0, \infty)$, $q\in (0, \pi/2)$, $r\in (0, \log(2+\sqrt{3}))$ and $s\in (0, \pi/3)$ be the parameters such that $1/\cosh(p)=\cos(q)=1-v^2$, $\cosh(r)=\sec(s)=1+v^2$. Then from (1.4)-(1.7) we have
\begin{equation}
\frac{\frac{1}{N_{AH}(a,b)}-\frac{1}{A(a,b)}}{\frac{1}{H(a,b)}-\frac{1}{A(a,b)}}
=\frac{\frac{2}{1+\frac{2p}{\sinh(2p)}}-1}{\frac{1}{1-v^{2}}-1}
=\frac{\frac{2\sinh(2p)}{2p+\sinh(2p)}-1}{\cosh(p)-1}=\varphi_{3}(p),
\end{equation}
\begin{equation}
\frac{\frac{1}{N_{HA}(a,b)}-\frac{1}{A(a,b)}}{\frac{1}{H(a,b)}-\frac{1}{A(a,b)}}
=\frac{\frac{2}{\cos(q)+\frac{q}{\sin(q)}}-1}{\frac{1}{1-v^{2}}-1}
=\frac{\frac{2\sin(q)-\sin(q)\cos(q)-q}{\sin(q)\cos(q)+q}}{\frac{1-\cos(q)}{\cos(q)}}
=1+\varphi_{4}(q),
\end{equation}
\begin{equation}
\frac{\frac{1}{N_{CA}(a,b)}-\frac{1}{C(a,b)}}{\frac{1}{A(a,b)}-\frac{1}{C(a,b)}}
=\frac{\frac{2}{\cosh(r)+\frac{r}{\sinh(r)}}-\frac{1}{1+v^{2}}}{1-\frac{1}{1+v^{2}}}=\varphi_{3}(r)
\end{equation}
and
\begin{equation}
\frac{\frac{1}{N_{AC}(a,b)}-\frac{1}{C(a,b)}}{\frac{1}{A(a,b)}-\frac{1}{C(a,b)}}=1+\varphi_{4}(s),
\end{equation}
where the functions $\varphi_{3}$ and $\varphi_{4}$ are defined as in Lemmas 2.6 and 2.8, respectively.

Note that
\begin{equation}
\varphi_{3}[\log(2+\sqrt{3})]=\frac{2\sqrt{3}-\log(2+\sqrt{3})}{2\sqrt{3}+\log(2+\sqrt{3})}
\end{equation}
and
\begin{equation}
\varphi_{4}\left(\frac{\pi}{3}\right)=-\frac{8\pi-6\sqrt{3}}{4\pi+3\sqrt{3}}.
\end{equation}

Therefore, Theorem 1.3 follows easily from (3.7)-(3.12) together with Lemmas
2.6 and 2.8. $\Box$

\medskip
\noindent{\bf Competing interests}

\noindent{The authors declare that they have no competing interests.}

\medskip
\noindent{\bf Authors' contributions}

\noindent{All authors contributed equally to the writing of this paper. All authors read and approved the final manuscript.}

\medskip
\noindent{\bf Acknowledgements}

\noindent{This research was supported by the
Natural Science Foundation of China under Grants 61374086, 11371125 and 11171307, and
the Natural Science Foundation of Zhejiang Province under Grant LY13A010004.}

\end{document}